\documentclass[a4paper,oneside,11pt]{article}

\usepackage[utf8]{inputenc}
\usepackage[T1]{fontenc}
\usepackage{textcomp}
\usepackage[english]{babel}
\usepackage[autolanguage]{numprint}
\usepackage{hyperref}
\usepackage{graphicx}
\usepackage{verbatim}
\usepackage{color}

\usepackage{amsmath,amssymb,amsthm}

\renewcommand\[{\begin{equation}}\renewcommand\]{\end{equation}}
\renewcommand\epsilon\varepsilon
\renewcommand\phi\varphi
\renewcommand\geq\geqslant
\renewcommand\leq\leqslant
\renewcommand\ln\log

\newcommand\Id{\operatorname{Id}}

\newcommand\ab\allowbreak

\theoremstyle{definition}

\theoremstyle{plain}

\newtheorem{thm}{Theorem}[section]

\newtheorem{lem}[thm]{Lemma}

\newtheorem{rmq}{Remark}[section]

\DeclareMathOperator{\Hess}{Hess}

\newcommand{\Z}{\mathbb{Z}}
\newcommand{\R}{\mathbb{R}}

\newcommand{\E}{\mathbb{E}}

\begin{document}

\title{HWI inequalities in discrete spaces via couplings}
\date{\today}

\author{Thomas A. Courtade\thanks{University of California, Berkeley, United States, Department of Electrical Engineering and Computer Sciences} 
 \hspace{1mm} and Max Fathi\thanks{Université Paris Cité and Sorbonne Université, CNRS, Laboratoire Jacques-Louis Lions and Laboratoire de Probabilit\'es, Statistique et Mod\'elisation, F-75013 Paris, France\\
and DMA, École normale supérieure, Université PSL, CNRS, 75005 Paris, France \\
and Institut Universitaire de France\\}}
\maketitle

\begin{abstract}
HWI inequalities are interpolation inequalities relating entropy, Fisher information and optimal transport distances. We adapt an argument of Y. Wu for proving the Gaussian HWI inequality via a coupling argument to the discrete setting, establishing new interpolation inequalities for the discrete hypercube and the discrete torus. In particular, we obtain an improvement of the modified logarithmic Sobolev inequality for the discrete hypercube of Bobkov and Tetali. 
\end{abstract}

\section{Introduction}

The HWI inequality, originally proved by Otto and Villani \cite{OV}, is an interpolation inequality relating entropy, Fisher information and $L^2$ transport distances. It plays a role in the synthetic theory of Ricci curvature bounds on metric spaces, and has found applications to concentration of measure, statistical physics and geometry. Alternative proofs have been given in \cite{Wu, GLRT20, KST}, and an improved dimensional version was derived in \cite{BGGK}. The most common point of view on this inequality is to view it as a consequence of the convexity of the entropy functional along certain families of interpolations, which can be interpreted as geodesics for a formal Riemanian structure on the space of probability measures. The main goal of this work is to explain how the proof of \cite{Wu} (see also \cite{BGRS14, RS13}), which does not rely on this convex viewpoint, can be adapted to the discrete setting. This leads to new discrete interpolation inequalities, different from those obtained by adapting the ideas of \cite{OV} (as was done in \cite{EM12}). 

In the Euclidean setting, the HWI inequality takes the following form: given a reference probability measure $d\mu = e^{-V}dx$ on $\mathbb{R}^d$ such that $\Hess V \geq K\Id$ for some $K \in \R$, we have for all other probability measures $\nu$ on $\mathbb{R}^d$
\begin{equation} \label{hwi_class}
H(\nu|\mu) \leq W_2(\nu, \mu)\sqrt{I(\nu|\mu)} - \frac{K}{2}W_2(\nu, \mu)^2
\end{equation}
where $W_2$ stands for the $L^2$ Wasserstein (or Monge-Kantorovitch) distance, $H$ for the relative entropy functional and $I$ for the relative Fisher information; formal definitions of each will be given later. When $K > 0$, the HWI inequality \eqref{hwi_class} implies a logarithmic Sobolev inequality, one of the main functional inequalities used for establishing concentration of measure estimates, as well as bounds on the trend to equilibrium for stochastic dynamics. Beyond their relationship with other functional inequalities, HWI inequalities have found some direct applications in statistical physics \cite{Fat13, HM14}. They hold in a more general setting of weighted manifold satisfying a Ricci curvature bound. A dimensional reinforcement of the Gaussian HWI inequality for even, strongly log-concave arguments was derived in \cite[Theorem 5.4]{AR23}

In the discrete setting, several families of HWI inequalities have been proposed \cite{EM12, GRST2014, Kra16}, as consequences of various proposed definitions of Ricci curvature bounds adapted to discrete spaces. For each of them, the approach consists in defining a family of interpolating curves in the space of probability measures, and proving that the entropy is uniformly (semi)-convex along these curves. These curves are interpreted as geodesic curves, and the distance used in place of the Wasserstein distance is the associated geodesic distance, while the Fisher information is the squared norm of the gradient of the entropy with respect to that metric structure. Such convexity properties are usually known as Ricci curvature bounds, in analogy with the Lott-Sturm-Villani synthetic notion of Ricci curvature bounds for Riemannian manifolds. All of these approaches have been shown to work for the simple graphs we shall consider here, and some of them have been shown to work for more sophisticated examples \cite{FM15, EHMT, Kra16}. 

Our starting point for this work is an alternative proof of the HWI inequality for the Gaussian space, due to Yihong Wu \cite{Wu}, which we shall describe in some detail in Section 2. What is interesting is that, while it uses some ingredients related to curvature bounds in the continuous setting (namely, a decay rate for the Fisher information along a stochastic dynamic), strictly speaking it does not require such a bound, and instead also relies on some very explicit coupling arguments to bound the entropy along the dynamic by the Wasserstein distance between the initial data and the equilibrium measure.  In particular, there is no need to introduce some family of geodesic interpolations. We shall mimic this proof in the discrete setting for several examples, and obtain new HWI-type inequalities, that are \emph{different} from those obtained using discrete curvature arguments. The distances involved will be simple variations on the  $L^1$ and $L^2$ Wasserstein distances, rather than the more sophisticated variational distances appearing for example in \cite{EM12, Kra16}. In the case of the hypercube, the inequality we obtain improves on the modified logarithmic Sobolev inequality of Bobkov and Tetali \cite{BT, PS16}.

\begin{rmq}
While this work was being finalized, Altschuler and Chewi published a preprint \cite{AC23} which also leverages couplings and convexity of entropy to prove reverse transport inequalities.  Their viewpoint offers some additional  flexibility by iterating short-time estimates, combined with regularity assumptions on the transition rates.  Unlike the discrete settings emphasized here, their work focuses on diffusion processes in the continuous setting, leveraging tools from stochastic calculus (e.g., Girsanov's theorem). 
\end{rmq}

\section{Yihong Wu's proof of the Gaussian HWI inequality}

The goal of this section is to present Y.~Wu's proof of the HWI inequality for the Gaussian measure, and extract the main arguments we will need to replicate in the discrete setting. This proof was included in \cite{BGRS14} and \cite{RS13}. We also note that a related approach was used in \cite{ORW12}, and extended in \cite{Shao} to establish an HWI inequality for the Wiener measure, using an infinite-dimensional Harnack inequality instead of a coupling argument. 

\subsection{The Gaussian HWI inequality}

Let $\gamma$ denote the standard Gaussian measure on $\mathbb{R}^d$, and let $\nu$ be another probability measure with  density $d\nu = \rho d\gamma$.  The relative entropy, Fisher information, and $W_2$ distance are defined respectively by
\begin{align*}
H(\nu|\gamma) &= \int \rho \log \rho d\gamma\\
I(\nu | \gamma) &= \int \frac{|\nabla \rho|^2}{\rho}d\gamma\\
W_2(\nu,\gamma)&= \inf_{X\sim \nu, Z\sim \gamma} \left( \E |X-Z|^2\right)^{1/2},
\end{align*}
where the infimum in the definition of $W_2$ is over all couplings of $X,Z$ with laws $X\sim \nu$ and $Z\sim \gamma$. If $\nu \not\ll \gamma$ we set $H(\mu|\gamma) = +\infty$, and we similarly set $I(\mu|\gamma) = +\infty$ if $\nu \not\ll \gamma$   or if $\rho$ is not weakly differentiable. It will be convenient to adopt the usual abuse of notation and write $H(X|Z)$ in place of  $H(\nu|\gamma)$ when $X\sim \nu$ and $Z\sim \gamma$, and similarly write $I(X|Z)\equiv I(\nu|\gamma)$ and $W_2(X,Z) \equiv W_2(\nu,\gamma)$.  In this notation, the Gaussian HWI inequality is 
\begin{align}
 H(X|Z) \leq W_2(X, Z)\sqrt{I(X|Z)} - \frac{1}{2}W_2(X,Z)^2. \label{eq:GaussianHWI}
\end{align}


To start the proof, we introduce the Ornstein--Uhlenbeck dynamic
$$dX_t = -X_tdt + \sqrt{2}dB_t, ~~t\geq 0$$
for a standard Brownian motion $(B_t)_{t\geq 0}$, 
which is a reversible diffusion process whose invariant distribution is the standard Gaussian measure. For initial data $X_0 \overset{\mathrm{law}}{=} X\sim \nu$, we also have the identity in law, or Mehler formula, 
\begin{equation} \label{mehler_gauss}
X_t \overset{\mathrm{law}}{=} e^{-t}X + \sqrt{1-e^{-2t}}Z,
\end{equation}
where $Z\sim \gamma$ is independent of $X$.  A direct computation shows that if $X$ and $Y$ are two Gaussians, with respective mean $x$ and $y$ and same variance $\sigma^2$, then
\begin{align}
H\left( X |  Y\right) = \frac{|x - y|^2}{2\sigma^2}.\label{relEntropyGaussians}
\end{align}
Relative entropy is jointly convex in its arguments, so for  $Z,Z'\overset{\mathrm{i.i.d.}}{\sim}\gamma$, the above implies
\begin{align*}
H(X_t|Z) &= H\left(e^{-t}X + \sqrt{1-e^{-2t}}Z | e^{-t}Z' + \sqrt{1-e^{-2t}}Z\right)\\
&\leq \frac{e^{-2t}}{2(1-e^{-2t})} \E[|X-Z'|^2]
\end{align*}
for every coupling of $X$ and $Z'$.  Taking the infimum over all such couplings, we get the following reverse transport-entropy inequality along the Ornstein--Uhlenbeck flow:  
\begin{equation} \label{rev_TE_OU}
H(X_t | Z) \leq \frac{e^{-2t}}{2(1 - e^{-2t})}W_2(X, Z)^2, ~~t\geq 0.
\end{equation}
This estimate also holds for non-Gaussian reference measures with curvature bounded from below \cite{BGL14}, with a different proof, but this is not our goal here. 

The classical entropy production formula (i.e., the de Bruijn identity) states 
\begin{equation*}
H(X|Z) = H(X_t|Z) + \int_0^t{I(X_s|Z)ds}, 
\end{equation*}
and since the Fisher information decays exponentially fast along the Ornstein-Uhlenbeck flow \cite{BE85, BGL14}, that is $I(X_t|Z) \leq e^{-2t}I(X|Z)$, we conclude 
\begin{equation*}
H(X|Z)\leq H(X_t|Z) + \frac{1-e^{-2t}}{2}I(X|Z) \leq \frac{e^{-2t}}{2(1 - e^{-2t})}W_2(X, Z)^2  + \frac{1-e^{-2t}}{2}I(X|Z).
\end{equation*}
Optimizing over $t\geq0$ produces the Gaussian HWI inequality \eqref{eq:GaussianHWI}. 

\subsection{Roadmap for a  scheme to prove HWI inequalities} \label{sec:roadmap}

The above argument relies on the fact that our reference measure $\gamma$ in the HWI inequality is the invariant measure of a certain Markov process, and that the Fisher information arrises as the derivative of the entropy along the process. Adopting this viewpoint, we can extract the two main ingredients of the proof we just described: 

\begin{enumerate}
\item A decay rate, or at least a bound, for the Fisher information along the flow generated by the Markov chain. 

\item A coupling of two trajectories of the Markov chain with different starting points, and such that the relative entropy of one marginal of the joint distribution at time $t$ with respect to the other can be estimated. 
\end{enumerate}

As we will illustrate in subsequent sections, when these two ingredients are available for some Markov process, they can be   combined to obtain a HWI inequality for the invariant measure.

The first element is one of the typical outcomes of a bound on entropic Ricci curvature, and so is available for any Markov chain that satisfies such a bound. But Fisher information decay is a strictly weaker property, and is sometimes known in cases where no lower bound on the Ricci curvature is known, or known but with worse constants. At a practical level, checking it requires checking convexity of the entropy along the dynamic flow, while Ricci curvature bounds require checking convexity along a much larger class of curves. See for example \cite{CDP} and \cite{FM15} to see how the two notions differ on concrete examples.  This part of the argument only involves Fisher information and entropy, it has no effect on what distance arises in the final inequality. In the examples below, we always use exponential decay rates for Fisher information, but in principle other rates could be used. 

The second element is at the core of the reverse transport-entropy inequality \eqref{rev_TE_OU}, and is established here using the Mehler formula \eqref{mehler_gauss}. Using couplings is a now-standard method for studying long-time behavior of Markov processes, see for example \cite{BD97}. It is in this step that the Wasserstein distance appears. 

In the discrete situations we shall consider next, we  demonstrate how these two elements can be used to mimic the proof of the Gaussian HWI inequality.  This approach leads to  different HWI inequalities than those obtained when only using curvature arguments. The downside of the method presented here is that an exact representation of transition probabilities (e.g., as given by the Mehler formula), is often unavailable, and thus the argument seems less widely applicable than curvature arguments. It may be that Harnack-type inequalities could be used to cover examples where we do not have explicit couplings, yet still lead to different HWI inequalities than those obtained via curvature arguments. 

\section{The hypercube}

Consider the simple random walk $(X_t)_{t\geq 0}$ on the hypercube $\{0,1\}^N$, where with rate $N$ we choose a coordinate uniformly at random and flip it with probability $1/2$. The invariant measure is the product Bernoulli measure with parameter $1/2$, which we shall denote by $\mu$. Given   initial data $X_0\equiv X \sim \nu$, we shall denote by $\nu_t$ the distribution after running this dynamic for some time $t$ (i.e., $X_t \sim \nu_t$). 

If we differentiate the entropy along this dynamic, we have
$$\frac{d}{dt}H(\nu_t|\mu) = -\frac{1}{2}\sum_{x\in \{0,1\}^N} (\rho_t(x^i) - \rho(x))(\log \rho_t(x^i) - \log \rho_t(x))\mu(x) =: I(\nu_t|\mu)$$
where $\rho_t$ is the density of $\nu_t$ with respect to $\mu$, and $x^i$ is the configuration of $\{0,1\}^N$ obtained from $x$ by flipping the $i$-th coordinate. The right-hand side is the discrete (modified) Fisher information, also known as entropy production. Note that under this scaling, the Fisher information $I$ is additive on the dimension. 

The main result of this section is the following HWI inequality for the discrete hypercube.  To state it, we define the $W_1$ distance 
$$
W_1(\nu,\mu) = \inf_{X\sim \nu, Y\sim \mu}\E[d(X,Y)],
$$
where $d$ is  the Hamming (graph) distance on the hypercube, and the infimum is over all couplings of $X\sim \nu$ and $Y\sim \mu$.
\begin{thm} Let $\mu$ be the uniform probability measure on $\{0,1\}^N$.  For any other probability measure $\nu$ on the same space satisfying $I(\nu|\mu) \geq 4W_1(\nu, \mu)$, we have
$$H(\nu|\mu) \leq 2\sqrt{W_1(\nu, \mu) I(\nu|\mu)} - 2W_1(\nu, \mu).$$
\end{thm}

Note that in the admissible regime $I\geq 4W_1$, this improves on the usual modified log-Sobolev inequality $H \leq \frac{1}{2}I$ of \cite{BT}, via Young's inequality. 

The transport-information inequality for the discrete hypercube states that $W_1^2 \leq N I$ (see for example \cite{FS18}), which is not enough to ensure $I \geq 4W_1$ is always true, but does imply it when $I/N$ is not too large.

\begin{rmq}
In this discrete setting, there are several notions of discrete Fisher information, and the associated functional inequalities are not equivalent in general, see for example \cite{BT}. 
\end{rmq}

\begin{proof}

It is known that $I$ is exponentially decaying: $I(\nu_t|\mu) \leq e^{-2t}I(\nu|\mu)$, which gives the first of the two ingredients discussed in Section \ref{sec:roadmap}. This can be derived for example as a consequence of the Ricci curvature bounds for the hypercube obtained in \cite{EM12} (but can also be computed directly). Hence 
$$H(\nu|\mu) \leq H(\nu_t|\mu) + \frac{1}{2}(1 - e^{-2t})I(\nu|\mu).$$
The usual proof of the modified logarithmic Sobolev inequality by the Bakry-Emery method consists at this point of letting $t$ go to infinity, to obtain $H \leq \frac{1}{2}I$.  Our  HWI inequality shall be obtained by taking a better choice for $t$, after bounding the entropy $H(\nu_t|\mu)$. 

We  now bound $H(\nu_t|\mu)$ using a coupling argument, giving the second ingredient in the discussion of Section \ref{sec:roadmap}. To do this, we couple two trajectories of the random walk starting from different deterministic initial conditions $X_0 = x$ and $Y_0 = y$ by having their same coordinates change at the same times, with same outcome. If we then look at the random variables $X_t$ and $X_t$ obtained by running this coupling, their coordinates that matched initially still match at later times, while coordinates that did not match at the beginning are the same with probability $1- e^{-2t}$. Moreover, their laws are products on the coordinates. Hence $H(X_t | X_t) \leq d(x,y)\varphi(t)$, where $\varphi(t)$ is the entropy of a Bernoulli random variable with parameter $p= (1 - e^{-2t})/2$ with respect to another one with parameter $q = (1+ e^{-2t})/2$. This quantity can be computed, and is
$$\varphi(t) = e^{-2t}\log\left(\frac{1+ e^{-2t}}{1 - e^{-2t}}\right).$$

Now, for any coupling of $X = X_0 \sim \nu$ and $Y = Y_0 \sim \mu$, we  have 
$$H(\nu_t |\mu) \leq   \varphi(t) \mathbb{E}[d(X,Y)]$$
and hence, minimizing over couplings, 
$$H(\nu_t|\mu) \leq \varphi(t)W_1(\nu, \mu).$$
Combining the above estimates, we get
$$H(\nu|\mu) \leq \varphi(t)W_1(\nu, \mu)+ \frac{1}{2}(1 - e^{-2t})I(\nu|\mu).$$

We can then use the bound $\varphi(t) \leq 2/(1 - e^{-2t}) - 2$ and take $t$ such that $(1 - e^{-2t}) = 2\sqrt{W_1}/I$ which is possible when $I \geq 4W_1$, to complete the proof. 
\end{proof}
\begin{rmq}
Note the difference with the Gaussian case, where it is a squared distance that plays a role. The appearance of $W_1$  makes sense, though, since it is $W_1$, and not $W_1^2$, that is additive on product measures. 
\end{rmq}
\section{The discrete torus}
In this section, we prove a new HWI inequality for the discrete torus.  We show that it implies the (known) HWI inequality on the circle as a limiting case.

\subsection{HWI inequality on discrete torus}

We now consider the situation where the reference measure $\mu$ is the uniform measure on the discrete hypercube $\Z/(N\Z)$, viewed as the invariant measure of the simple random walk. The relative entropy of another probability measure $\nu$ is then  
$$H(\nu|\mu) = \sum_{x \in \Z/(N\Z)} \nu(x)(\log \nu(x) + \log N),$$
and the Fisher information is
$$I(\nu|\mu) := \sum_{x \in \Z/(N\Z)} (\log \nu(x+1) - \log \nu(x))(\nu(x+1) - \nu(x)),$$
which is indeed the dissipation of entropy along the flow of the simple random walk.  

The HWI inequality for the uniform measure on the discrete torus of length $N$ will involve the following transport cost
$$W_c(\nu, \mu)^2 := \underset{X\sim \nu, Y\sim \mu}{\inf} \hspace{1mm} \mathbb{E}[d(X,Y) + d(X,Y)^2],$$
where once again the infimum is taken over all possible couplings of $X\sim \nu$ and $Y\sim \mu$, and $d$ denotes the graph distance. Note that
$$\max(W_1, W_2^2) \leq W^2_c \leq \min(2W_2^2, (N+1)W_1).$$

We shall obtain the following HWI inequality: 

\begin{thm}\label{HWIdiscreteTorus}
Let $\mu$ be the uniform probability measure on $\Z/(N\Z)$.  For any other probability measure $\nu$ on the same space, we have the HWI inequality
$$H(\nu|\mu) \leq \sqrt{2}W_c(\nu,\mu)\sqrt{I(\nu|\mu)}.$$
\end{thm}

\begin{rmq}
The HWI inequality for the continuous torus obtained when adopting the viewpoint of \cite{OV} is $H \leq W_2\sqrt{I}$. One can recover this inequality from the discrete one above by rescaling. It is also possible to prove the continuous HWI inequality directly with the method we use here by coupling two Brownian motions on the torus. 
\end{rmq}

\begin{proof}
Let $X_t$ and $Y_t$ be two simple random walks on the discrete torus starting from positions at distance $d$. We realize them starting from simple random walks $\tilde{X}_t$ and $\tilde{Y}_t$ on the integers also starting at distance $d$, and setting $X_t = \tilde{X}_t \text{ mod } N$ (resp. $Y_t = \tilde{Y}_t \text{ mod } N$). Starting with the data processing inequality for relative entropy, we have
$$H(X_t |Y_t) \leq H(\tilde{X}_t | \tilde{Y}_t) = e^{-t}\underset{n \in \Z}{\sum} \hspace{1mm} I_n(t)\log\left(\frac{I_n(t)}{I_{n-d}(t)}\right)$$
where $I_n(t)$ is the so-called modified Bessel function of the first kind \cite{AS}, which is related to the transition probabilities of the simple random walks via $\mathbb{P}(X_t = n | X_0 = 0) = e^{-t}I_n(t)$. From Lemma \ref{cor_bessel} below, we obtain
$$H(X_t | Y_t) \leq \frac{d + d^2}{2t}.$$

For random walks with general initial data $X_0 = X$ and $Y_0 = Y$, the same coupling argument as before yields
$$H(X_t | Y_t) \leq \frac{W_c(X, Y)^2}{2t}$$
and hence, since the Fisher information is non-increasing along the flow \cite{EM12}, 
$$H(\nu|\mu) \leq \frac{W_c(\nu, \mu)^2}{2t} + t I(\nu|\mu).$$
Optimizing in $t$ gives the result. 
\end{proof}

We now prove the technical estimate on transition probabilities of the simple random walk we used above: 

\begin{lem} \label{cor_bessel}
Let $M$ be a symmetric, unimodal integer-valued random variable. Then 
$$\mathbb{E}\left[\log \frac{I_M(t)}{I_{M-d}(t)}\right] \leq \frac{d + d^2}{2t}.$$
\end{lem}

For a symmetric probability measure $\mathbb{P}$, unimodality means that $m \longrightarrow \mathbb{P}(m)$ is non-increasing in $|m|$. 

\begin{proof}
This lemma is a consequence of the following estimate on modified Bessel functions: for integers $n, d$ satisfying $n \geq d/2 \geq 0$, we have for all $t > 0$
\begin{equation} \label{bnd_ratio_bessel}
\log \frac{I_n(t)}{I_{n-d}(t)} \geq \frac{(1 +d)(d-2n)}{2t}.
\end{equation}
Let us take this bound as given for now, and show how it implies Lemma \ref{cor_bessel}. For simplicity, we assume that $d$ is even, so that $d/2$ is an integer. Let $h(n) := \frac{(1 +d)(d-2n)}{2t} - \log \frac{I_n(t)}{I_{n-d}(t)}$. Since it is the difference of two functions that are odd about the point $d/2$, it also is. Moreover, as a consequence of \eqref{bnd_ratio_bessel}, $h(n) \leq 0$ for $n \geq d/2$, and by antisymmetry  $h(n) \geq 0$ for $n \leq d/2$. 

We then have
\begin{align*}
\mathbb{E}[h(M)] &= \sum_{k \geq 1} \mathbb{P}(M = d/2 -k)h(d/2-k) + \mathbb{P}(M = d/2 +k)h(d/2+k) \\
&= \sum_{k \geq 1} (\mathbb{P}(M = d/2 -k) - \mathbb{P}(M = d/2 +k)h(d/2-k) 
\end{align*}
When $k \geq 0$, $h(d/2-k) \geq 0$, and moreover $\mathbb{P}(M = d/2 -k) \geq \mathbb{P}(M = d/2 +k)$ by the assumptions on the distribution of $M$. Hence $\mathbb{E}[h(M)] \geq 0$, and therefore
$$\mathbb{E}\left[\log \frac{I_M(t)}{I_{M-d}(t)}\right] = \mathbb{E}\left[\frac{(1 +d)(d-2M)}{2t}\right] - \mathbb{E}[h(M)] \leq \frac{d + d^2}{2t}.$$
The case where $d$ is odd follows the same chain of arguments. 

Let us now prove \eqref{bnd_ratio_bessel}. In \cite[p.241]{Am74}, it is shown that
\begin{equation}
\sqrt{1 + \left(\frac{n+1}{t}\right)^2} - \frac{n+1}{t} \leq \frac{I_{n+1}(t)}{I_n(t)}.
\end{equation}
We define 
$$f(t,x) := \log\left(\sqrt{1 + \left(\frac{x+1}{t}\right)^2} - \frac{x+1}{t}\right) - \left(\frac{x^2}{tc(x)} - \frac{(x+1)^2}{tc(x+1)}\right).$$
with $c(x) = 2(1 - (x+1)^{-1})$. Direct computations show that for all $t, x \geq 0$ 
$$\partial_x f(t,x) \geq 0; \hspace{3mm} \partial_t f(t,0) \leq 0$$
and therefore, for any $t \geq 0$ and $x \geq 0$
$$f(t, x) \geq \lim_{t \to \infty} f(t,0) = 0.$$
In particular, 
$$\log \frac{I_{n+1}(t)}{I_n(t)} \geq \frac{n^2}{c(n)t} - \frac{(n+1)^2}{c(n+1)t}.$$
For $n \geq d$, it immediately follows that
$$\log \frac{I_{n}(t)}{I_{n-d}(t)} \geq \frac{(n-d)^2}{c(n-d)t} - \frac{n^2}{c(n)t} = \frac{d(d-1-2n)}{2t} \geq \frac{(1+d)(d-2n)}{2t}.$$
Now, if $d/2 \leq n < d$, then $n \geq d-n > 0$ and 
$$\log \frac{I_{n}(t)}{I_{n-d}(t)} = \log \frac{I_{n}(t)}{I_{d-n}(t)} \geq \frac{(1+d)(d-2n)}{2t}.$$
\end{proof}

\vskip2ex

\textbf{\underline{Acknowledgments}} : This work was started while the second author was visiting the MSRI in Berkeley during the semester program \emph{Geometric Functional Analysis} in the Fall 2017. Support from the France-Berkeley Fund, NSF-CCF 1750430,  ANR-11-LABX-0040-CIMI within the program ANR-11-IDEX-0002-02 is acknowledged. This work was partially supported by Projects EFI (ANR-17-CE40-0030) and CONVIVIALITY (ANR-23-CE40-0003) of the French National Research Agency, and the FSMP under the Invited Professor program.

\end{document}